\documentclass[11pt]{article}
\usepackage{graphicx}
\usepackage{amsmath,amsthm,amsfonts,amssymb}
\newtheorem{teo}{Theorem}[section]
\newtheorem{coro}[teo]{Corollary}
\newtheorem{lema}[teo]{Lemma}
\newtheorem{prop}[teo]{Proposition}
\newtheorem{obs}[teo]{Remark}%[section]
\newtheorem{ejem}[teo]{Example}

\def\qed{\begin{flushright} \QED \end{flushright}}

\newcommand\ZZ{{\mathbb{Z}}}

\def\C{{\mathcal C}}

\def\F{{\mathcal F}}

\def\L{{\mathcal L}}
\def\U{{\mathcal U}}

\def\Hom{\mathop{\sf Hom}\nolimits}

\def\qed{\hfill \mbox{$\square$}\bigskip}

\def\Ext{\mathop{\sf Ext}\nolimits}

\def\Ker{\mathop{\rm Ker}\nolimits}

\def\Im{\mathop{\rm Im}\nolimits}

\def\mod{\mathop{\sf mod}\nolimits}

\begin{document}

\sf  \title{$\mathsf{N}$-complexes as functors, amplitude cohomology and fusion rules}
\author{Claude Cibils, Andrea Solotar and Robert Wisbauer
\thanks{\footnotesize This work has been supported by the projects PICT 08280 (ANPCyT), UBACYTX169, PIP-CONICET 5099 and
 the German Academic Exchange Service (DAAD). The second author is a research member of
CONICET (Argentina) and a Regular Associate of ICTP Associate Scheme.}}

\date{}

\maketitle

\begin{abstract}
We consider $\mathsf{N}$-complexes as functors over an appropriate
linear category in order to show first that the Krull-Schmidt
Theorem holds, then to prove that amplitude cohomology (called
generalized cohomology by M. Dubois-Violette) only vanishes on
injective functors providing a well defined functor on the stable
category. For left truncated $\mathsf{N}$-complexes, we show that
amplitude cohomology discriminates the isomorphism class up to a
projective functor summand. Moreover amplitude cohomology of
positive $\mathsf{N}$-complexes is proved to be isomorphic to an
$\Ext$ functor of an indecomposable $\mathsf{N}$-complex inside the
abelian functor category. Finally we show that for the monoidal
structure of $\mathsf{N}$-complexes a Clebsch-Gordan formula holds,
in other words the fusion rules for $\mathsf{N}$-complexes can be
determined.

\end{abstract}



%%%%%%%%%%%%%%%%%%%%%%%%%%%%%%%%%%%%%%%%%%%%%%%%%%%%%%%%%%%%%%%%%%%%%%%%%%%%%%%%%%%%%%%%%%
\section{Introduction}

Let $\mathsf{N}$ be a positive integer and let $k$ be a field. In
this paper we will consider $\mathsf{N}$-complexes of vector spaces
as linear functors (or modules) over a $k$-category, see the
definitions at the beginning of Section \ref{categories}.

Recall first that a usual $k$-algebra is deduced from any finite
object $k$-category through the direct sum of its vector spaces of
morphisms. Modules over this algebra are precisely $k$-functors from
the starting category, with values in the category of $k$-vector
spaces. Consequently if the starting category has an infinite number
of objects, linear functors with values in vector spaces are called
modules over the category, as much as modules over an algebra are
appropriate algebra morphisms.

An $\mathsf{N}$-complex as considered by M. Kapranov in \cite{ka} is
a $\ZZ$-graded vector space equipped with linear maps $d$ of degree
$1$ verifying $d^\mathsf{N}=0$. The amplitude (or generalized)
cohomology are the vector spaces $\Ker d^a/\Im d^{\mathsf{N}-a}$ for
each amplitude $a$ between $1$ and $\mathsf{N}-1$. Note that we use
the terminology \emph{amplitude cohomology} in order to give a
graphic idea of this theory and in order to clearly distinguish it
from classical cohomology theories.

M. Dubois-Violette has shown in \cite{du} a key result, namely that
for $\mathsf{N}$-complexes arising from cosimplicial modules through
the choice of an element $q\in k$ such that
$1+q+\cdots+q^{\mathsf{N}-1}=0$, amplitude cohomology can be
computed using the classical cohomology provided the truncated sums
$1+q+\cdots+q^{n}$ are invertible for $1\leq n\leq \mathsf{N}-1$. As
a consequence he obtains in a unified way that Hochschild cohomology
at roots of unity or in non-zero characteristic is zero or
isomorphic to classical Hochschild cohomology (see also \cite{kawa})
and the result proven in 1947 by Spanier \cite{sp}, namely that
Mayer \cite{ma} amplitude cohomology can be computed by means of
classical simplicial cohomology.

Note that $\mathsf{N}$-complexes are useful for different
approaches, as Yang-Mills algebras \cite{codu}, Young symmetry of
tensor fields \cite{duhe,duhe2} as well as for studying homogeneous
algebras and Koszul properties, see
\cite{beduwa,bema,dupo,masa,masa2} or for analysing cyclic homology
at roots of unity \cite{wa}. A comprehensive description of the use of
$\mathsf{N}$-complexes in this various settings is given in
the course by M.
Dubois-Violette at the Institut Henri Poincar\'{e}, \cite{du.ihp}.

We first make clear an obvious fact, namely that an
$\mathsf{N}$-complex is a module over a specific $k$-category
presented as a free $k$-category modulo the $\mathsf{N}$-truncation
ideal. This way we obtain a Krull-Schmidt theorem for
$\mathsf{N}$-complexes. The list of indecomposables is well-known,
in particular projective and injective $\mathsf{N}$-complexes
coincide. This fact enables us to enlarge Kapranov's aciclicity
Theorem  in terms of injectives. More precisely, for each amplitude
 $a$ verifying $1\le a \le \mathsf{N}-1$
a classic $2$-complex is associated to each $\mathsf{N}$-complex. We
prove first in this paper that an $\mathsf{N}$-complex is acyclic
for a given amplitude if and only if the $\mathsf{N}$-complex is
projective (injective), which in turn is equivalent to aciclicity
for any amplitude.

In \cite{duke,du} a basic result is obtained for amplitude
cohomology for $\mathsf{N}\geq 3$ which has no counterpart in the
classical situation $\mathsf{N}=2$, namely hexagons raising from
amplitude cohomologies are exact. This gap between the classical and
the new theory is confirmed by a result we obtain in this paper:
amplitude cohomology does not discriminate arbitrary
$\mathsf{N}$-complexes without projective summands for $\mathsf{N}
\ge 3$, despite the fact that for $\mathsf{N}=2$ it is well known
that usual cohomology is a complete invariant up to a projective
direct summand. Nevertheless we prove that left truncated
$\mathsf{N}$-complexes sharing the same amplitude cohomology are
isomorphic up to a projective (or equivalently injective) direct
summand.

We also prove that amplitude cohomology for \emph{positive}
$\mathsf{N}$-complexes coincides with an $\Ext$ functor in the
category of $\mathsf{N}$-complexes. In other words, for each given
amplitude there exists an indecomposable module such that the
amplitude cohomologies of a positive $\mathsf{N}$-complex are
actually extensions of a particular degree between the
indecomposable and the given positive $\mathsf{N}$-complex. We use
the characterisation of $\Ext$ functors and the description of
injective positive $\mathsf{N}$-complexes. In this process the fact
that for positive $\mathsf N$-complexes, projectives no longer
coincide with injectives requires special care.

We underline the fact that various indecomposable modules are used in
order to show that amplitude cohomology of positive
$\mathsf{N}$-complexes is an $\Ext$ functor. This variability makes
the result compatible with the non classical exact hexagons
\cite{duke,du} of amplitude cohomologies quoted above.

M. Dubois-Violette has studied in \cite{du.bariloche} (see Appendix
A) the monoidal structure of $\mathsf{N}$-complexes in terms of the
coproduct of the Taft algebra, see also \cite{du.ihp}. J. Bichon in
\cite{bi} has studied the monoidal structure of
$\mathsf{N}$-complexes, considering them as comodules, see also the
work by R. Boltje \cite{bo} and A. Tikaradze \cite{ti}. We recall in
this paper that the $k$-category we consider is the universal cover
of the Taft Hopf algebra $\U^+_q(sl_2)$. As such, there exists a
tensor product of modules (i.e. $\mathsf{N}$-complexes) for each
non-trivial $\mathsf N^{\mathsf{th}}$-root of unity (see also
\cite{bo,ci}). Using Gunnlaugsdottir's axiomatisation of
Clebsch-Gordan's formula \cite{gu} and amplitude cohomology we show
that this formula is valid for $\mathsf{N}$-complexes, determining
this way the corresponding fusion rules.

%%%%%%%%%%%%%%%%%%%%%%%%%%%%%%%%%%%%%%%%%%%%%%%%%%%%%%%%%%%%%%%%%%%%%%%%%%%%%%%%%%%%%%%%%%%
\section{$\mathsf{N}$-complexes and categories}
\label{categories}

Let $\C$ be a small category over a field $k$. The set of objects is denoted $\C_{0}$.
Given $x,y$ in $\C_0$, the $k$-vector space of morphisms from $x$ to $y$ in $\C$ is
denoted ${}_{y}\C_{x}$ . Recall that  composition of morphisms is $k$-bilinear. In this
way, each ${}_{x}\C_{x}$ is a $k$-algebra and each ${}_{y}\C_{x}$ is a
${}_{y}\C_{y}$-${}_{x}\C_{x}$ -- bimodule.

For instance let $\Lambda$ be a $k$-algebra and let $E$ be a complete finite system of
orthogonal idempotents in $\Lambda$, that is $\sum_{e\in E} e =1$, $ef=fe=0$ if $f\neq e$
and $e^2=e$, for all $e,f \in E$. The associated category $\C_{\Lambda,E}$ has set of
objects $E$ and morphisms ${}_{f}\left(\C_{\Lambda,E}\right)_{e}=f\Lambda e$. Conversely
any finite object set category $\C$ provides an associative algebra $\Lambda$ through the
matrix construction. Both procedures are mutually inverse.

In this context linear functors $F: \C_{\Lambda,E} \to Mod_k$
coincide with left $\Lambda$-modules. Consequently for any arbitrary
linear category $\C$, left modules are defined as $k$-functors $F:
\C \to Mod_k$. In other words, a left $\C$-module is a set of
$k$-vector spaces $\{{}_xM\}_{x\in \C_0}$ equipped with "left
oriented" actions that is, linear maps
\[{}_y\C_x \otimes_k {}_xM \to {}_yM\]
verifying the usual associativity constraint.

Notice that right modules are similar, they are given by a collection of $k$-vector
spaces $\{M_x\}_{x\in \C_0}$  and "right oriented"' actions. From now on a module will
mean a left module.

Free $k$-categories are defined as follows: let $E$ be an arbitrary set and let $V=
\{{}_yV_x\}_{x,y \in E} $ be a set of $k$-vector spaces. The free category $\F_E(V)$  has
set of objects $E$ and set of morphisms from $x$ to $y$ the direct sum of tensor products
of vector spaces relying $x$ to $y$:

\[ {}_y(\F_E(V))_x =
\bigoplus_{n\ge 0}\  \ \bigoplus_{x_1,\dots,x_n \in E}({}_yV_{x_n}\otimes \dots \otimes
{}_{x_2}V_{x_1}\otimes {}_{x_1}V_x ) \]

For instance, let $E=\ZZ$ and let ${}_{i+1}V_i =k$ while ${}_jV_i=0$
otherwise. This data can be presented by the double infinite quiver
having $\ZZ$ as set of vertices and an arrow from $i$ to $i+1$ for
each $i\in \ZZ$. The corresponding free category $\L$ has one
dimensional vector space morphisms from $i$ to $j$ if and only if
$i\le j$, namely
\[{}_j\L_i= {}_jV_{j-1}\otimes \dots \otimes {}_{i+2}V_{i+1}\otimes {}_{i+1}V_i .\]

Otherwise ${}_j\L_i=0$.

A module over $\L$ is precisely a graded vector space $\{{}_iM\}_{i\in \ZZ}$ together
with linear maps $d_i:{}_iM \to {}_{i+1}M$. This fact makes use of the evident universal
property characterizing free linear categories.

On the other hand we recall from \cite{ka} the definition of an $\mathsf{N}$-complex: it
consists of a graded vector space  $\{{}_iM\}_{i\in \ZZ}$ and linear maps $d_i:{}_iM \to
{}_{i+1}M$ verifying that $d_{i+\mathsf{N}}\circ \dots \circ d_i=0$ for each $i\in \ZZ$.

In order to view an $\mathsf{N}$-complex as a module over a $k$-linear category we have
to consider a quotient of $\L$. Recall that an ideal $I$ of a $k$-category $\C$ is a
collection of sub-vector spaces ${}_yI_x$ of each morphism space ${}_y\C_x$, such that
the image of the composition map ${}_z\C_y \otimes {}_yI_x$ is contained in ${}_zI_x$ and
${}_yI_x \otimes {}_x\C_u$ is contained in ${}_yI_u$ for each choice of objects. Quotient
$k$-categories exist in the same way that algebra quotients exist.

Returning to the free category $\L$, consider the truncation ideal $I_\mathsf{N}$ given
by the entire ${}_j\L_i$ in case $j\ge i+\mathsf{N}$ and $0$ otherwise. Then
$\L_\mathsf{N}:=\L/I_\mathsf{N}$ has one dimensional morphisms from $i$ to $j$ if and
only if $i\le j \le i+\mathsf{N}-1$.

Clearly $\mathsf{N}$-complexes coincide with $\L_\mathsf{N}$-modules. We have obtained
the following

\begin{teo}
The categories  of $\mathsf{N}$-complexes and of $\L_\mathsf{N}$-modules are isomorphic.
\end{teo}

An important point is that $\L_\mathsf{N}$ is a locally bounded category, which means
that the direct sum of morphism spaces starting (or ending) at each given object is
finite dimensional. More precisely:

$$ \forall x_0, y_0 \in (\L_\mathsf{N})_0, \ \ dim_k \left[\bigoplus_{y \in \ZZ}{}_y(\L_\mathsf{N})_{x_0}\right]
 = \mathsf{N}= dim_k \left[\bigoplus_{x \in
\ZZ}{}_{y_0}(\L_\mathsf{N})_{x}\right].$$

It is known that for locally bounded categories Krull-Schmidt theorem holds, for instance
see the work by C. S\'aenz \cite{sa}. We infer that each $\mathsf{N}$-complex of finite
dimensional vector spaces is isomorphic to a direct sum of indecomposable ones in an
essentially unique way, meaning that given two decompositions, the multiplicities of
isomorphic indecomposable $\mathsf{N}$-complexes coincide.

Moreover, indecomposable $\mathsf{N}$-complexes are well known, they correspond to "short
segments" in the quiver: the complete list of indecomposable modules is given by
$\{M_i^l\}_{i \in \ZZ, 0\le l \le \mathsf{N}-1}$ where $i$ denotes the beginning of the
module, $i+l$ its end and $l$ its length. More precisely, ${}_i(M_i^l)={}_{i+1}(M_i^l)=
\dots = {}_{i+l}(M_i^l)=k$ while ${}_j(M_i^l)= 0$ for other indices $j$. The action of
$d_i, d_{i+1}, \dots , d_{i+l-1}$ is the identity and $d_j$ acts as zero if the index $j$
is different. The corresponding $\mathsf{N}$-complex is concentrated in the segment
$[i,i+l]$.

Note that the simple $\mathsf{N}$-complexes are $\{M_i^0\}_{i \in \ZZ}$ and that each
$M_i^l$ is uniserial, which means that $M_i^l$ has a unique filtration
\[ 0\subset M_{i+l}^0 \subset \dots \subset M_{i+2}^{l-2} \subset M_{i+1}^{l-1}
\subset M_i^l \]
such that each submodule is maximal in the following one.

Summarizing the preceding discussion, we have the following

\begin{prop}
Let $M$ be an $\mathsf{N}$-complex of finite dimensional vector spaces. Then
\[ M \simeq \ \ \bigoplus_{i \in \ZZ, \ 0\le l\le \mathsf{N}-1}n_i^lM_i^l \]
for a unique finite set of positive integers $n_i^l$.
\end{prop}

Indecomposable projective and injective $\L_\mathsf{N}$-modules are also well known, we
now recall them briefly. Note from \cite{fre} that projective functors are direct sums of
representable functors. Clearly ${}_{-}(\L_\mathsf{N})_i=M_i^{\mathsf{N}-1}$.

In order to study injectives notice first that for a locally bounded $k$-category, right
and left modules are in duality: the dual of a left module is a right module which has
the dual vector spaces at each object, the right actions are obtained by dualising the
left actions. Projectives and injectives correspond under this duality. Right projective
modules are direct sums of ${}_{i}(\L_\mathsf{N})_-$ as above, clearly
$({}_{i}(\L_\mathsf{N})_-)^* \simeq M_i^{\mathsf{N}-1}$.

This way we have provided the main steps of the proof of the following

\begin{prop}
Let $M_i^l$ be an indecomposable $\mathsf{N}$-complex, $i\in \ZZ$ and $l\le
\mathsf{N}-1$. Then $M_i^l$ is projective if and only if $l=\mathsf{N}-1$, which in turn
is equivalent for $M_i^l$ to be injective.
\end{prop}

\begin{coro}
Let $M= \bigoplus_{i \in \ZZ, 0\le l\le \mathsf{N}-1}n_i^lM_i^l$ be an
$\mathsf{N}$-complex. Then $M$ is projective if and only if $n_i^l=0$ for $l\le
\mathsf{N}-2$, which in turn is equivalent for $M$ to be injective.
\end{coro}

%%%%%%%%%%%%%%%%%%%%%%%%%%%%%%%%%%%%%%%%%%%%%%%%%%%%%%%%%%%%%%%%%%%%%%%%%%%%%%%%%%%%%%%%%%%
\newpage
\section{Amplitude cohomology}

Let $M$ be an $\mathsf{N}$-complex. For each amplitude $a$ between $1$ and
$\mathsf{N}-1$, at each object $i$ we have $\Im d^{\mathsf{N}-a}\subseteq \Ker d^a$. More
precisely we define as in \cite{ka}
\[(AH)_a^i(M):= \Ker (d_{i+a-1}\circ \dots \circ d_i)/\Im(
d_{i-1}\circ \dots \circ d_{i-\mathsf{N}+a}) \] and we call this
bi-graded vector space the amplitude cohomology of the
$\mathsf{N}$-complex. As remarked in the Introduction, M.
Dubois-Violette in \cite{du} has shown the depth of this theory, he
calls it generalised cohomology.

As a fundamental example we compute amplitude cohomology for
indecomposable $\mathsf{N}$-complexes $M_i^l$. In the following
picture the amplitude is to be read vertically while the
 degree of the cohomology is to be read horizontally. A black dot
means one dimensional cohomology, while an empty dot stands for zero
cohomology.

\hskip1cm
\includegraphics[width=1\textwidth,height=0.6\textheight]{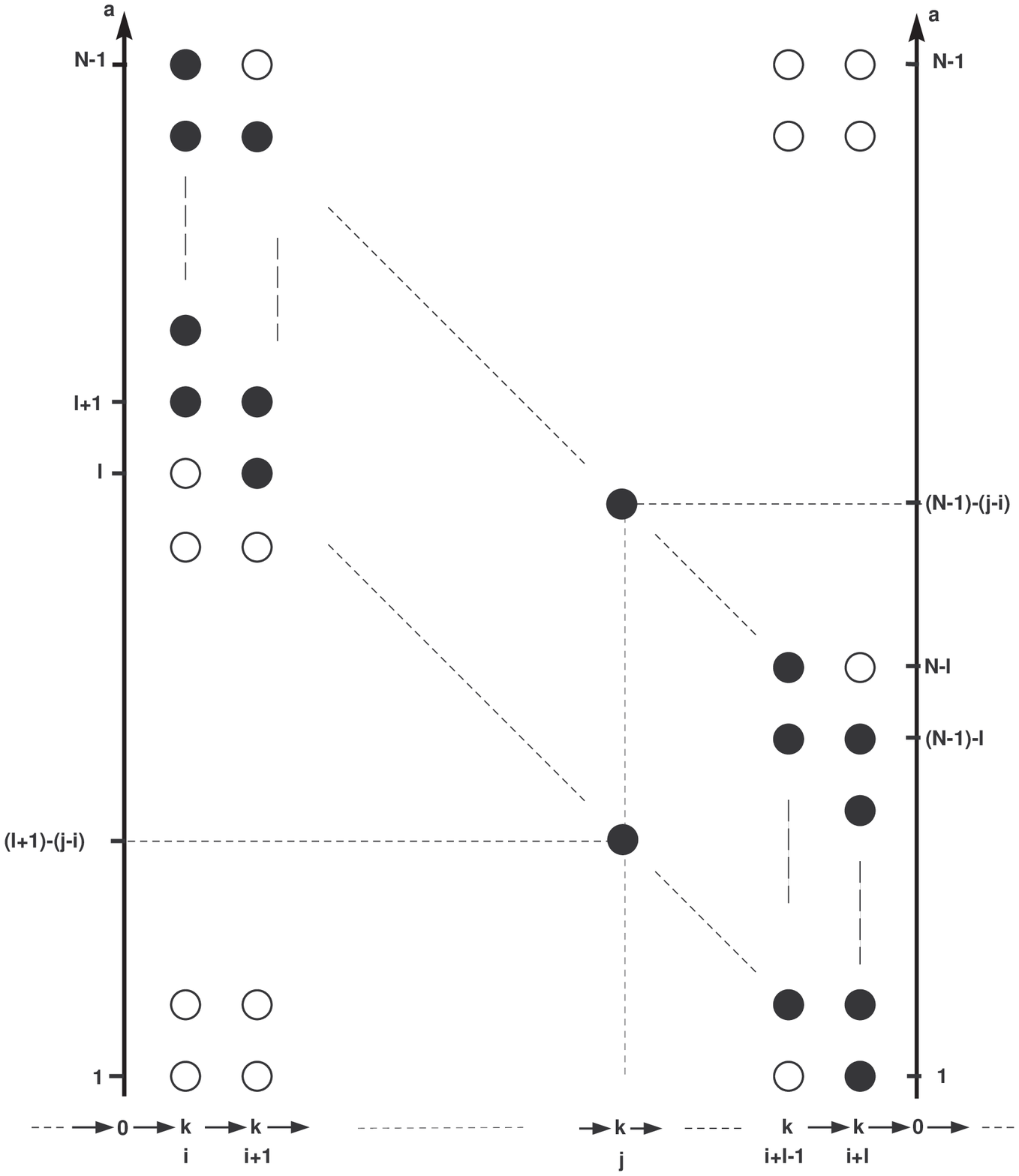}
From this easy computation we notice that for a non-projective
(equivalently non-injective) indecomposable module $M_i^l$ ($0\le
l\le \mathsf{N}-2$) and any amplitude $a$ there exists a degree $j$
such that $(AH)^j_a(M_i^l) \neq 0$. Concerning projective or
injective indecomposable modules $M_i^{\mathsf{N}-1}$ we notice that
$(AH)^j_a(M_i^{\mathsf{N}-1})=0$ for any degree $j$ and any
amplitude $a$. These facts are summarized as follows:

\begin{prop}
Let $M$ be an indecomposable $\mathsf{N}$-complex. Then $M$ is projective (or
equivalently injective) in the category of $\mathsf{N}$-complexes if and only if its
amplitude cohomology vanishes at some amplitude $a$ which in turn is equivalent to its
vanishing at any amplitude.
\end{prop}

\begin{obs}
\sf From the very definition of amplitude cohomology one can check that for a fixed
amplitude $a$ we obtain a linear functor $(AH)_a^*$ from $\mod\L_\mathsf{N}$ to the
category of graded vector spaces.

Moreover $(AH)_a^*$ is additive, in particular:
\[ (AH)_a^*(M\oplus M') = (AH)_a^*(M) \oplus (AH)_a^*(M')\]
\end{obs}

This leads to the following result, which provides a larger frame to
the aciclicity result of M. Kapranov \cite{ka}. See also the short
proof of Kapranov's aciclicity result by M. Dubois-Violette in Lemma
3 of \cite{du} obtained as a direct consequence of a key result of
this paper, namely the exactitude of amplitude cohomology hexagons.

\begin{teo} Let $M$ be an $\mathsf{N}$-complex of finite dimensional vector
spaces. Then $(AH)_a^*(M)=0$ for some $a$ if and only if $M$ is
projective (or equivalently injective). Moreover, in this case
$(AH)_a^*(M)=0$ for any amplitude
 $a\in [1, N-1]$.
\end{teo}

In order to understand the preceding result in a more conceptual framework we will
consider the stable category of $\mathsf{N}$-complexes, $\underline{\mod\L_\mathsf{N}}$.
More precisely, let $I$ be the ideal of $\mod\L_\mathsf{N}$ consisting of morphisms which
factor through a projective $\mathsf{N}$-complex. The quotient category
$\mod\L_\mathsf{N}/I$ is denoted $\underline{\mod\L_\mathsf{N}}$. Clearly all projectives
become isomorphic to zero in $\underline{\mod\L_\mathsf{N}}$. Of course this construction
is well known and applies for any module category. We have in fact proven the following

\begin{teo}
For any amplitude $a$ there is a well-defined functor
\[(AH)_a^*: \underline{\mod\L_\mathsf{N}} \to gr(k)\]
where $gr(k)$ is the category of graded $k$-vector spaces.
\end{teo}

Our next purpose is to investigate how far amplitude cohomology distinguishes
$\mathsf{N}$-complexes. First we recall that in the classical case ($\mathsf{N}=2$),
cohomology is a complete invariant of the stable category.

\begin{prop}
Let $M$ and $M'$ be $2$-complexes of finite dimensional vector spaces
without projective direct summands.
If $H^*(M) \simeq H^*(M')$, then $M \simeq M'$.
\end{prop}
\noindent\textbf{Proof.} Indecomposable $2$-complexes are either simple or projective. We
assume that $M$ has no projective direct summands, this is equivalent for $M$ to be
semisimple, in other words $M$ is a graded vector space with zero differentials.
Consequently $H^i(M)={}_iM$ for all $i$. \qed

The following example shows that the favorable situation for $\mathsf{N}=2$ is no longer
valid for $\mathsf{N}\ge 3$.

\begin{ejem} \sf
Consider $M$ the $3$-complex which is the direct sum of all simple modules, in other
words, ${}_iM=k$ and $d_i=0$. Then for any  degree $i$ we have
\[  (AH)_1^i(M)=k \hbox{ and }(AH)_2^i(M)=k.  \]
Let $M'$ be the direct sum of all the length one indecomposable $3$-complexes,
$$M'=\bigoplus_{i\in \ZZ} M^1_i$$
Recall that the amplitude cohomology of $M_i^1$ is given by
\[ (AH)_2^i(M_i^1)=k \hbox{ and }(AH)_1^{i+1}(M_i^1)=k  \]
while all other amplitude cohomologies vanish. Summing up provides $(AH)_2^i(M')=k $ and
$(AH)_1^i(M')=k $, for all $i$. However it is clear that $M$ and $M'$ are not isomorphic.
Notice that both $M$ and $M'$ are free of projective direct summands.
\end{ejem}

As quoted in the introduction the preceding example confirms that
amplitude cohomology is a theory with different behaviour than the
classical one. This fact has been previously noticed by M.
Dubois-Violette in \cite{du}, for instance when dealing with  non
classical exact hexagons of amplitude cohomologies.

At the opposite, we will obtain in the following that for either
left or right truncated $\mathsf{N}$-complexes amplitude cohomology
is a complete invariant up to projectives. More precisely, let $M$
be an $\mathsf{N}$-complex which is zero at small enough objects,
namely ${}_iM=0$ for $i\le b$, for some $b$ which may depend on $M$.
Of course this is equivalent to the fact that for the Krull-Schmidt
decomposition
\[ M=\bigoplus_{i\in \ZZ}\bigoplus_{l=0}^{\mathsf{N}-1}n_i^lM_i^l \]
there exists a minimal $i_0$, in the sense that
$n_i^l=0$ if $i< i_0$ and $n_{i_0}^l
\neq 0$ for some $l$.

\begin{prop}
Let $M$ be a non-projective $\mathsf{N}$-complex which is zero at small enough objects.
Let $l_0$ be the smallest length of an indecomposable factor of $M$ starting at the
minimal starting object $i_0$. Then $(AH)_a^{i}(M)=0$ for all $i\le i_0-1$ and
$(AH)_a^{i_0}(M)=0$ for $a\le l_0$. Moreover $dim_k(AH)_{l_0+1}^{i_0}(M)=n_{i_0}^{l_0}$.
\end{prop}
\noindent\textbf{Proof.} The fundamental computation we made of amplitude cohomology for
indecomposable $\mathsf{N}$-complexes shows the following: the smallest degree affording
non vanishing amplitude cohomology provides the starting vertex of an indecomposable non
projective module. Moreover, at this degree the smallest value of the amplitude affording
non zero cohomology is $l+1$, where $l$ is the length of the indecomposable.

In other words amplitude cohomology determines the multiplicity of the \emph{smallest}
indecomposable direct summand of a left-truncated $\mathsf{N}$-complex. Of course
\emph{smallest} concerns the lexicographical order between indecomposables, namely $M_i^l
\le M_j^r$ in case $i<j$ or in case $i=j$ and $l\le r$. \qed

\begin{teo}
\label{algo} Let $M$ be an $\mathsf{N}$-complex which is zero at
small enough objects and which does not have projective direct
summands. The dimensions of its amplitude cohomology determine the
multiplicities of each indecomposable direct summand.
\end{teo}
\noindent\textbf{Proof.}
The Proposition above shows that the multiplicity of the smallest indecomposable direct summand is
determined by the amplitude cohomology
(essentially this multiplicity is provided by the smallest non-zero
amplitude cohomology, where amplitude cohomology
is also ordered by lexicographical order).

We factor out this smallest direct summand $X$ from $M$ and we notice that the
multiplicities of other indecomposable factors remain unchanged. Moreover, factoring out
the amplitude cohomology of $X$ provides the amplitude cohomology of the new module. It's
smallest indecomposable summand comes strictly after $X$ in the lexicographical order.
Through this inductive procedure, multiplicities of indecomposable summands can be
determined completely. In other words: if two left-truncated $\mathsf{N}$-complexes of
finite dimensional vector spaces share the same amplitude cohomology, then the
multiplicities of their indecomposable direct factors coincide for each couple $(i,n)$.
\qed

\begin{obs}
\sf Clearly the above Theorem is also true for $\mathsf{N}$-complexes which are zero for
large enough objects, that is right-truncated $\mathsf{N}$-complexes.

\end{obs}

%%%%%%%%%%%%%%%%%%%%%%%%%%%%%%%%%
\section{ Amplitude cohomology is $\Ext$}

An $\mathsf{N}$-complex $M$ is called positive in case ${}_iM=0$ for $i\le -1$. In this
section we will prove that amplitude cohomology of positive $\mathsf{N}$-complexes of
finite dimensional vector spaces coincides with an $\Ext$ functor in this category.

First we provide a description of injective positive $\mathsf{N}$-complexes as modules.
Notice that positive $\mathsf{N}$-complexes are functors on the full subcategory
$\L_\mathsf{N}^{\geqslant 0}$ of $\L_\mathsf{N}$ provided by the positive integer
objects. Alternatively, $\L_\mathsf{N}^{\geqslant 0}$ is the quotient of the free
$k$-category generated by the quiver having positive integer vertices and an arrow from
$i$ to $i+1$ for each object, by the truncation ideal given by morphisms of length
greater than $\mathsf{N}$.

\begin{teo}
The complete list up to isomorphism of injective positive indecomposable
$\mathsf{N}$-complexes is
\[ \{M_0^l\}_{l=0,\dots,\mathsf{N}-1} \sqcup \{M_i^{\mathsf{N}-1}\}_{i\ge 1} \]
\end{teo}
\noindent\textbf{Proof.} As we stated before, injective modules are duals of projective
right modules. The indecomposable ones are representable functors
${}_{i_0}(\L_\mathsf{N}^{\geqslant 0})_{-}$, for $i_0\ge 0$.

Clearly for each $i_0$ we have $({}_{i_0}{\L_\mathsf{N}^{\geqslant 0}}_{-})^* =
M_0^{i_0}$ if $i_0\le \mathsf{N}-1$ while $({}_{i_0}{\L_\mathsf{N}^{\geqslant 0}}_{-})^*
= M_{i_0-(\mathsf{N}-1)}^{\mathsf{N}-1}$ otherwise.

In order to show that amplitude cohomology is an instance of an $\Ext$, we need to have
functors sending short exact sequences of positive $\mathsf{N}$-complexes into long exact
sequences: this will enable to use the axiomatic characterization of $\Ext$. For this
purpose we recall the following standard consideration about $\mathsf{N}$-complexes which
provides several classical $2$-complexes associated to a given $\mathsf N$-complex, by
contraction. More precisely fix an integer $e$ as an initial condition and an amplitude
of contraction $a$  (which provides also a coamplitude of contraction $b=\mathsf{N}-a$).

The contraction $C_{e,a}M$ of an $\mathsf{N}$-complex is the
following $2$-complex, which has ${}_eM$ in degree $0$ and
alternating $a$-th and $b$-th composition differentials:

\[ \dots \rightarrow {}_{e-b}M \stackrel{d^{b}}{\rightarrow}
         {}_eM{} \stackrel{d^{a}}{\rightarrow} {}_{e+a}M \stackrel{d^{b}}{\rightarrow}
       {}_{e+\mathsf{N}}M \stackrel{d^{a}}{\rightarrow} \dots     \]

Of course usual cohomology of this complex provides amplitude cohomology:

\begin{lema}
In the above situation, $$H^{2i}(C_{e,a}M) = (AH)_a^{e+iN}(M) \mbox{
\ and\ \ } H^{2i+1}(C_{e,a}M) = (AH)_b^{e+iN+a}(M).$$

\end{lema}
Notice that in order to avoid repetitions and in order to set $H^0$
as the first positive degree amplitude cohomology, we must restrict
the range of the initial condition. More precisely, for a given
amplitude contraction $a$ the initial condition $e$ verifies $0\le e
< b$, where $b$ is the coamplitude verifying $a+b=\mathsf{N}$.
Indeed, if $e\ge b$, set $e'=e-b$ and $a'=b$. Then $b'=a$ and $0\le
e' < b'$.

\begin{obs}
\sf An exact sequence of $\mathsf{N}$-complexes provides an exact sequence of contracted
complexes at any initial condition $e$ and any amplitude $a$.
\end{obs}
We focus now on the functor $H^*(C_{e,a}-)$, which for simplicity we shall denote
$H^*_{e,a}$ from now on.  We already know that  $H^*_{e,a}$ sends a short exact sequence
of $\mathsf{N}$-complexes into a long exact sequence, since $H^*_{e,a}$ is usual
cohomology. Our next purpose is two-fold. First we assert that $H^*_{e,a}$ vanishes in
positive degrees when evaluated on injectives of the category of positive
$\mathsf{N}$-complexes. Then we will show that it is representable in degree 0.

\begin{prop}
In positive degrees we have: $$ H^*_{e,a}(M_0^l)=0 \mbox{ for } l\le
\mathsf{N}-1,\ \mbox{ and } \  H^*_{e,a}(M_i^{\mathsf{N}-1})=0
\hbox{ for } i\ge 1.$$

\end{prop}
\noindent\textbf{Proof.} Concerning indecomposable modules of length $\mathsf{N}-1$, they
are already injective in the entire category of $\mathsf{N}$-complexes. We have noticed
that all their amplitude cohomologies vanish.

Consider now $M_0^l$, with $l\le \mathsf{N}-1$. In non-zero even degree $2i$ the
amplitude cohomology to be considered is in degree $e+iN$, which is larger than $l$ since
$i\neq 0$ and
%%%%%%%%%%
$\mathsf{N}> l$. Hence $H^{2i}_{e,a}(M_0^l)=0 $.

In odd degree $2i+1$ the amplitude cohomology to be considered is in degree $e+iN+a$. As
before, in case $i\neq 0$ this degree is larger than $l$, then $H^{2i+1}_{e,a}(M_0^l)=0 $
for $i\neq 0$. It remains to consider the case $i=0$, namely $H^{1}_{e,a}(M_0^l)=
(AH)_{\mathsf{N}-a}^{e+a}(M_0^l)$. From the picture we have drawn for amplitude
cohomology in the previous section, we infer that in degree $e+a$ the cohomology is not
zero only for amplitudes inside the closed interval $[l+1-(e+a),\ \ \mathsf{N}-1-(e+a)]$.
We are concerned by the amplitude $\mathsf{N}-a$ which is larger than $\mathsf{N}-a-e-1$,
hence $H^{1}_{e,a}(M_0^l)= 0$. \qed

\begin{prop}
Let $a \in [1,\mathsf{N}-1]$ be an amplitude and let $e \in [0,\mathsf{N}-1-a]$ be an
initial condition. Then $H^{0}_{e,a}(-)= (AH)_a^e(-)$ is a representable  functor given
by the indecomposable $\mathsf{\mathsf{N}}$-complex $M_e^{a-1}$. More precisely,

\[ (AH)_a^e(X) = \Hom_{\L_\mathsf{N}^{\geqslant 0}}(M_e^{a-1}, X). \]
\end{prop}
\noindent\textbf{Proof.} We will verify this formula for an arbitrary indecomposable
positive $\mathsf{N}$-complex $X=M_i^l$. The morphism spaces between indecomposable
$\mathsf{N}$-complexes are easy to determine using diagrams through the defining quiver
of $\L_\mathsf{N}^{\geqslant 0}$. Non-zero morphisms from an indecomposable $M$ to an
indecomposable $M'$ exist if and only if $M$ starts during $M'$ and $M$ ends together
with or after $M'$. Then we have:

\[ \Hom_{\L_\mathsf{N}^{\geqslant 0}}(M^{a-1}_{e}, M_i^l) = \begin{cases}
                                                          k & \text{if $e \in [i,i+l]$ and $e+a-1\ge i+l$ }
                                                                    \\
                                                          0 & \text{otherwise}
                                               \end{cases}
   \]

Considering amplitude cohomology and the fundamental computation
we have made, we first notice that $(AH)_a^e(M_i^l)$ has a chance to be non-zero
only when the degree $e$ belongs to the indecomposable, namely
$e \in [i,i+l]$. This situation already coincides with the first condition
for non-vanishing of $\Hom$. Next, for a given $e$ as before, the precise
conditions that the amplitude $a$ must verify in order to obtain $k$ as
amplitude cohomology is
\[ (l+1)-(e-i) \le a \le (\mathsf{N}-1)-(e-i). \]

The second inequality holds since the initial condition $e$ belongs to
$[0,\mathsf{N}-1-a]$ and $i\ge0$. The first inequality is precisely $ e+a-1 \ge i+l$.
\qed

As we wrote before it is well known (see for instance \cite{macl}) that a functor
sending naturally short exact sequences into long exact sequences,
vanishing on injectives and being representable in degree $0$ is isomorphic
to the corresponding $\Ext$ functor.
Then we have the following:

\begin{teo}
Let $\L_\mathsf{N}^{\geqslant 0}$ be the category of positive $\mathsf{N}$-complexes of
finite dimensional vector spaces and let $AH_a^j(M)$ be the amplitude cohomology of an
$\mathsf{N}$-module $M$ with amplitude $a$ in degree $j$. Let $b=\mathsf{N}-a$ be the
coamplitude.

Let $j=q\mathsf{N}+e$ be the euclidean division with $0\le e \le \mathsf{N}-1$.

Then for $e < b$ we have:

\[ AH_a^j(M)= \Ext_{\L_\mathsf{N}^{\geqslant 0}}^{2q}(M_{e}^{a-1}, M).  \]

and for $e\ge b$ we have:

\[ AH_a^j(M)= \Ext_{\L_\mathsf{N}^{\geqslant 0}}^{2q+1}(M_{e-b}^{b-1}, M).  \]

\end{teo}

%%%%%%%%%%%%%%%%%%%%%%%%%%%%%%

\section{Monoidal structure and Clebsch-Gordan formula}

The $k$-category $\L_\mathsf{N}$ is the universal cover of the associative algebra
$\U_q^+(sl_2)$ where $q$ is a non-trivial $\mathsf{N}$-th root of unity, see \cite{ci} and also \cite{ciro}.
More precisely, let $C=<t>$ be the infinite cyclic group and let $C$ act on
$(\L_\mathsf{N})_0= \ZZ$ by $t.i=i+N$. This is a free action on the objects while the
action on morphisms is obtained by translation: namely the action of $t$ on the generator
of ${}_{i+1}V_{i}$ is the generator of ${}_{i+1+N}V_{i+\mathsf{N}}$.

Since the action of $C$ is free on the objects, the categorical quotient exists, see for
instance \cite{cire}. The category $\L_\mathsf{N}/C$ has set of objects $\ZZ/\mathsf{N}$.
This category $\L_\mathsf{N}/C$ has a finite number of objects, hence we may consider its
matrix algebra $a(\L_\mathsf{\mathsf{N}}/C)$ obtained as the direct sum of all its
morphism spaces equipped with matrix multiplication.  In other words,
$a(\L_\mathsf{N}/C)$ is the path algebra of the crown quiver having $\ZZ/\mathsf{N}$ as
set of vertices and an arrow form $\bar i$ to $\bar i+1$ for each $\bar i \in
\ZZ/\mathsf{N}$, truncated by the two-sided ideal of paths of length
%%%%%%%%%%%%%%%
greater or equal to
$\mathsf{N}$.

As described in \cite{ci} this truncated path algebra bears a comultiplication, an
antipode and a counit providing a Hopf algebra isomorphic to the Taft algebra, also known
as the positive part $\U_q^+(sl_2)$ of the  quantum group $\U_q(sl_2)$. The monoidal
structure obtained for the $\U_q^+(sl_2)$-modules can be lifted to
$\L_\mathsf{N}$-modules providing the monoidal structure on $\mathsf{N}$-complexes
introduced by M. Kapranov \cite{ka} and studied by J. Bichon \cite{bi} and A. Tikaradze \cite{ti}.

We recall the formula: let $M$ and $M'$ be $\mathsf{N}$-complexes. Then $M\otimes M'$ is
the $\mathsf{N}$-complex given by

\[ {}_i(M\otimes M') = \bigoplus_{j+r=i}({}_jM\otimes {}_rM') \]

and

\[ d_i(m_j \otimes m'_r) = m_j \otimes d_rm'_r + q^rd_jm_j\otimes m'_r. \]

Notice that in general ${}_i(M\otimes M')$ is not finite dimensional.

\begin{prop}
Let $M$ and $M'$ be $\mathsf{N}$-complexes of finite dimensional vector spaces. Then
$M\otimes M'$ is a direct sum of indecomposable $\mathsf{N}$-complexes of finite
dimensional vector spaces, each indecomposable appearing a finite number of times.
\end{prop}
\noindent\textbf{Proof.}
Using Krull-Schmidt Theorem we have

\[ M= \bigoplus_{i \in \ZZ,\ 0\le l\le \mathsf{N}-1}n_i^lM_i^l \hbox{\hspace{0.1in} and \hspace{0.1in}  }
   M'= \bigoplus_{i \in \ZZ,\ 0\le l\le \mathsf{N}-1}n'^l_i M_i^l . \]

The tensor product $M_i^l\otimes M_j^r$ consists of a finite number of non-zero vector
spaces which are finite dimensional. It follows from the Clebsch-Gordan formula that we
prove below that for a given indecomposable $\mathsf{N}$-complex $M_l^u$, there is only a
finite number of couples of indecomposable modules sharing $M_l^u$ as an indecomposable
factor. Then each indecomposable appears a finite number of times in $M\otimes M'$. \qed

The following result is a Clebsch-Gordan formula for indecomposable
$\mathsf{N}$-complexes, see also the work by R. Boltje, chap. III \cite{bo}. The fusion
rules, \emph{i.e.} the positive coefficients arising from the decomposition of the tensor
product of two indecomposables, can be determined as follows.

\begin{teo}
Let $q$ be a non-trivial $\mathsf{N}^{th}$ root of unity and $M_i^u$ and let $M_j^v$ be
indecomposable $\mathsf{N}$-complexes. Then,\\
if $u+v\le \mathsf{N}-1$ we have
\[ M_i^u\otimes M_j^v = \bigoplus_{l=0}^{min(u,v)}M_{i+j+l}^{u+v-2l}, \]
if $u+v= e+\mathsf{N}-1$ with $e\ge 0$ we have
\[ M_i^u\otimes M_j^v = \bigoplus_{l=0}^{e}M_{i+j+l}^{\mathsf{N}-1}
 \oplus \bigoplus_{l=e+1}^{min(u,v)}M_{i+j+l}^{u+v-2l}.  \]

\end{teo}
\noindent\textbf{Proof.} Using Gunnlaugsdottir's axiomatization \cite{gu} p.188, it is
enough to prove the following:

\[ M_i^0\otimes M_j^0 = M_{i+j}^0 \]

\[ M_0^1\otimes M_j^u = M_{j}^{u+1} \oplus M_{j+1}^{u-1} \hbox{ for }u < \mathsf{N}-1 \]

 \[ M_0^1\otimes M_j^{\mathsf{N}-1} = M_{j}^{\mathsf{N}-1} \oplus M_{j+1}^{\mathsf{N}-1} . \]

The first fact is trivial. The second can be worked out using amplitude cohomology, which
characterizes truncated $\mathsf{N}$-complexes. Indeed the algorithm we have described in
\ref{algo} enables us first to determine the fusion rule for $M_0^1\otimes M_j^u $
($u<\mathsf{N}-1$), that is to determine the non-projective indecomposable direct
summands. More precisely, since $u < \mathsf{N}-1$, the smallest non vanishing amplitude
cohomology degree is $j$, with smallest amplitude $u+2$, providing $M_j^{u+1}$ as a
direct factor. The remaining amplitude cohomology corresponds to $M_{j+1}^{u-1}$. A
dimension computation shows that in this case there are no remaining projective summands.

On the converse, the third case is an example of vanishing cohomology. In fact, since
$M_j^{\mathsf{N}-1}$ is projective, it is known at the Hopf algebra level that $X \otimes
M_j^{\mathsf{N}-1}$ is projective. A direct dimension computation between projectives
shows that the formula holds. \qed

%%%%%%%%%%%%%%%%%%%%%%%%%%%%%%%%%%%%%%%%%%%%%%%%%%%%%%%%%%%%%%%%%%%%%%%%%%%%%%%%%%%%%%%%%%%

%%%%%%%%%%%%%%%%%%%%%%%%%%%%%%%%%%%%%%%%%%%%%%%%%%%%%%%%%%%%%%%%%%%%%%%%%%%%%%%%%%%%%%%%%%%

\footnotesize
\noindent C.C.:
\\ Institut de Math\'ematiques et de Mod\'{e}lisation de Montpellier I3M,\\
UMR 5149\\
Universit\'e de Montpellier 2,
\\F-34095 Montpellier cedex 5,
France.\\
{\tt Claude.Cibils@math.univ-montp2.fr}

\noindent A.S.:
\\Departamento de Matem\'atica,
 Facultad de Ciencias Exactas y Naturales,\\
 Universidad de Buenos Aires
\\Ciudad Universitaria, Pabell\'on 1\\
1428, Buenos Aires, Argentina. \\{\tt asolotar@dm.uba.ar}

\noindent R.W.:
\\ Mathematical Institute,
Heinrich-Heine-University\\
Universitaetsstrasse 1\\
D-40225 Duesseldorf,
Germany
\\{\tt  wisbauer@math.uni-duesseldorf.de}


\begin{thebibliography}{99}

\bibitem{beduwa} Berger, R.; Dubois-Violette, M.;
Wambst, M. Homogeneous algebras. J. Algebra \textbf{261} (2003), no. 1, 172–185.

\bibitem{bema} Berger, R.; Marconnet N. Koszul and Gorenstein Properties for Homogeneous
Algebras. Algebras and Representation Theory {\bf 9}, (2006),
67--97.



\bibitem {bi} Bichon, J.  {$\mathsf{N}$-complexes et alg\`ebres de Hopf}.
 C. R. Math. Acad. Sci. Paris  {\bf 337},  (2003),  441--444.

 \bibitem{bo} Boltje, R.,
      Kategorien von verallgemeinerten Komplexen
      und ihre Beschreibung durch Hopf Algebren,
      Diplomarbeit, Universit\"at M\"unchen (1985).

\bibitem {ci} Cibils, C. {A quiver quantum group }.
 Comm. Math. Phys. {\bf 157}, (1993),  459--477.

\bibitem {cire} Cibils, C.; Redondo M. J. {Cartan-Leray spectral
sequence for Galois coverings of categories}. J. Algebra {\bf 284}, (2005), 310--325.

\bibitem{ciro} Cibils, C.; Rosso, M. Hopf quivers.  J. Algebra  \textbf{254}  (2002),  241--251.


\bibitem{codu} Connes, A.; Dubois-Violette, M. Yang-Mills algebra.
Lett. Math. Phys.  \textbf{61}  (2002),  149--158.



\bibitem {du} Dubois-Violette, M. {$d^\mathsf{N}=0$: generalized homology}.
 $K$-Theory  {\bf 14},  (1998),   371--404.

\bibitem{du2} Dubois-Violette, M.
Generalized homologies for $d\sp \mathsf{N}=0$ and graded $q$-differential algebras.
Secondary calculus and cohomological physics (Moscow, 1997), 69--79,
Contemp. Math., \textbf{219}, Amer. Math. Soc., Providence, RI, 1998.

\bibitem{du.bariloche}
Dubois-Violette, M. Lectures on differentials, generalized
differentials and on some examples related to theoretical physics.
Quantum symmetries in theoretical physics and mathematics
(Bariloche, 2000),  59--94, Contemp. Math., \textbf{294}, Amer.
Math. Soc., Providence, RI, 2002.

\bibitem{du.ihp}
Dubois-Violette, M. R\'esum\'e et transparents du cours
"$\mathsf{N}$-COMPLEXES", for the semester "K-theory and
noncommutative geometry", Institut Henri Poincar\'e,
Paris, mars 2004.\\
http://qcd.th.u-psud.fr/page\_perso/MDV/articles/COURS\_IHP.pdf



\bibitem{duhe} Dubois-Violette, M.; Henneaux, M.
Generalized cohomology for irreducible tensor fields of mixed Young symmetry type.
Lett. Math. Phys. \textbf{49} (1999), 245--252.


\bibitem{duhe2} Dubois-Violette, M.; Henneaux, M.
Tensor fields of mixed Young symmetry type and $\mathsf{N}$-complexes.
Comm. Math. Phys. \textbf{226} (2002), 393--418.

\bibitem{duke} Dubois-Violette, M.; Kerner, R.: Universal q-differential calculus and q-analog of
homological algebra, Acta Math. Univ. Comenian. \textbf{65} (1996), 175–188.

\bibitem{dupo} Dubois-Violette, M.; Popov, T.
Homogeneous algebras, statistics and combinatorics.
Lett. Math. Phys. \textbf{61} (2002), 159--170.



\bibitem {fre} Freyd, P. \textit {Abelian categories}.
Harper and Row, New York, 1964.\\
\texttt{http://www.tac.mta.ca/tac/reprints/articles/3/tr3.pdf}

\bibitem {gu} Gunnlaugsdottir, E.  {Monoidal structure of
the category of $u\sp +\sb q$-modules}. Linear Algebra Appl. {\bf 365}, (2003), 183--199.

\bibitem {ka} Kapranov, M. On the q-analog of homological algebra. Preprint, Cornell University,
1991; q-alg/9611005

\bibitem {kawa} Kassel, C.; Wambst, M.  { Alg\`ebre
homologique des $\mathsf{N}$-complexes et homologie de Hochschild aux racines de
l'unit\'e }.
 Publ. Res. Inst. Math. Sci. {\bf 34}, (1998), 91--114.

\bibitem {macl} Mac Lane, S. \emph{Homology}, Springer-Verlag, New York, 1963.

\bibitem{ma} Mayer, W. A new homology theory. I, II.
Ann. of Math.  {\bf 43}, (1942). 370--380, 594--605.

\bibitem{masa} Mart\'{\i}nez-Villa, R.; Saor\'\i n, M. Koszul duality for $\mathsf N$-Koszul algebras.
Colloq. Math.  \textbf{103}  (2005), 155--168.

\bibitem{masa2} Mart\'{\i}nez-Villa, R.; Saor\'\i n, A duality theorem for generalized Koszul algebras\\
math.RA/0511157


\bibitem {sa} S\'aenz, C. \textit {Descomposici\'on en inescindibles
para m\'odulos sobre anillos y categor\'\i as asociadas}. Tesis
para obtener el t\'\i tulo de matem\'atico, UNAM, Mexico, 1988.

\bibitem{sp} Spanier, E.H. The Mayer homology theory, Bull. Amer. Math. Soc. 55 (1949),
102–-
112.

\bibitem{ti} Tikaradze, A. Homological constructions on $\mathsf{N}$-complexes.
J. Pure Appl. Algebra \textbf{176} (2002), 213–222.

\bibitem{wa} Wambst, M. Homologie cyclique et homologie simpliciale aux racines de l'unit\'e.
$K$-Theory  \textbf{23}  (2001), 377--397.

\end{thebibliography}
\end{document}